\documentclass{article}
\usepackage{amssymb,latexsym,amscd}

\newcommand{\se}[1]{{\section{#1}} {\setcounter{equation}{0}}}
\newtheorem{theorem}{Theorem}[section]
\newtheorem{lm}{Lemma}[section]
\newtheorem{prop}{Proposition}[section]

\newtheorem{co}{Corollary}[section]

\def\k{{K\"{a}hler }}

\begin{document}
\hbadness=10000
\title{{\bf Generalized special Lagrangian torus fibration for Calabi-Yau hypersurfaces in toric varieties II}}
\author{Wei-Dong Ruan\\
Department of Mathematics\\
University of Illinois at Chicago\\
Chicago, IL 60607\\}
\footnotetext{Partially supported by NSF Grant DMS-0104150.}
\maketitle
\begin{abstract}
In this paper we construct monodromy representing generalized special Lagrangian torus fibrations for Calabi-Yau hypersurfaces in toric varieties near the large complex limit.
\end{abstract}
\se{Introduction}
This paper is a sequel of \cite{sl1}. The aim of this series of papers is to construct generalized special Lagrangian torus fibrations for Calabi-Yau hypersurfaces in toric varieties near the large complex limit. In this paper we will use the same notations as in \cite{sl1} unless specified otherwise.\\
 
Let $(P_\Delta,\omega)$ be a toric variety whose moment map image (with respect to the toric \k form $\omega$) is the real convex polyhedron $\Delta \subset M_{\mathbb{R}}$. Also assume that the anti-canonical class of $P_\Delta$ is represented by an integral reflexive convex polyhedron $\Delta_0\subset M$ and the unique interior point of $\Delta_0$ is the origin of $M$. Integral points $m \in \Delta_0$ correspond to holomorphic toric sections $s_m$ of the anti-canonical bundle. For the unique interior point $m_o$ of $\Delta_0$, $s_{m_o}$ is the section of the anti-canonical bundle that vanishes to order 1 along each toric divisor of $P_\Delta$.\\

Let $\{w_m\}_{m\in \Delta_0}$ be a strictly convex function on $\Delta_0$ such that $w_m> 0$ for $m\in \Delta_0\setminus \{m_o\}$ and $w_{m_o}\ll 0$. Define

\[
\tilde{s}_t = s_{m_o} + ts,\ \ s = \sum_{m\in \Delta_0\setminus \{m_o\}} a_m s_m,\  {\rm where}\ |a_m| = \tau^{w_m},\ {\rm for}\ m\in \Delta_0\setminus \{m_o\}.
\]

(As in the case of Fermat type quintic, in general, $\Delta_0$ need not contain all the integral $m$ in the real polyhedron spanned by $\Delta_0$.) Let $X_t = \{\tilde{s}_t^{-1}(0)\}$. Then $\{X_t\}$ is a 1-parameter family of Calabi-Yau hypersurfaces in $P_{\Delta}$. $X_0 = \{s_{m_o}^{-1}(0)\}$ is the so-called large complex limit. $X_t$ is said to be near the large complex limit if $\tau$ and $t$ are small and $t \leq \tau^{-w_{m_o}}$.\\

Since $X_0$ is toric, the moment map induces the standard generalized special Lagrangian fibration $\hat{\pi}_0: X_0 \rightarrow \partial \Delta$ with respect to the toric holomorphic volume form. In \cite{lag1, lag2, lag3, toric}, we constructed Lagrangian torus fibration for $X_t$ when $X_t$ is near the large complex limit, using the Hamiltonian-gradient flow to deform this fibration $\hat{\pi}_0$ for $X_0$ symplectically to the desired Lagrangian fibration $\hat{\pi}_t: X_t \rightarrow \partial \Delta$ for such $X_t$. The (topological) singular set of the fibration map $\hat{\pi}_t$ is $C= X_t \cap {\rm Sing}(X_0)$, which is independent of $t$. The corresponding singular locus $\tilde{\Gamma} = \hat{\pi}_0(C)$ is also independent of $t$. When $X_t$ is near the large complex limit, $\tilde{\Gamma} \subset \partial \Delta$ exhibits amoeba structure that is a fattening of a graph $\Gamma \subset \partial \Delta$. It was conjectured in \cite{lag2} that the singular locus for the (generalized) special Lagrangian torus fibration should resemble the singular locus $\tilde{\Gamma}$ of the Lagrangian torus fibration.\\ 

\cite{sl1} initiated the program of using similar idea to construct generalized special Lagrangian fibration for $X_t$ by deformation from the standard fibration $\hat{\pi}_0$ of the large complex limit $X_0$. As pointed out in \cite{sl1}, due to the canonical nature of the generalized special Lagrangian fibration on a Calabi-Yau when we fix the \k metric, one can construct the generalized special Lagrangian fibration over different parts of the base $\partial \Delta$ separately and they will automatically match on the overlaps.\\ 

Recall from \cite{sl1} that a (generalized) special Lagrangian fibration of smooth torus over an open set $U\subset \partial \Delta$ is said to represent the monodromy if $\partial \Delta \setminus U$ is a fattening of $\tilde{\Gamma}$ that retracts to $\tilde{\Gamma}$. \\

There is a natural toric map from $P_\Delta$ to its anti-canonical model $P_{\Delta_0}$ that induces a map $\pi: \partial \Delta \rightarrow \partial \Delta_0$. Recall that $\Delta_0$ is an integral reflexive polyhedron. We will use $\Delta_0^{(k)}$ to denote the set of integral points in the $k$-skeleton of $\Delta_0$, where the coefficients of the polynomial defining the family of Calabi-Yau hypersurfaces are non-zero. Let $U \subset \partial \Delta$ be the union of the interior of top dimensional faces of $\partial \Delta$ and a tubular neighborhood of $\pi^{-1}(\Delta_0^{(n-1)})$ in $\partial \Delta$. (Here $n=\dim \partial \Delta$.) It is easy to observe that a fibration over $U$ represents monodromy from our discussion of the singular locus $\tilde{\Gamma}$ in \cite{lag3, toric}.\\

In \cite{sl1}, we constructed generalized special Lagrangian fibration for $X_t$ over two types of regions in $\partial \Delta$. The type one regions are the interior of top dimensional faces in $\partial \Delta$ (also see \cite{G}). The type two regions are small neighborhoods of vertices of $\partial \Delta$. Notice from \cite{lag3,toric} that these two types of regions are exactly dual to each other under the so-called symplectic SYZ mirror symmetry proved in \cite{toric} for Calabi-Yau hypersurfaces in toric varieties. This construction will be enough to produce monodromy representing generalized special Lagrangian fibration if $\Delta_0^{(n-1)}=\Delta_0^{(0)}$ and $\dim \pi^{-1}(\Delta_0^{(0)}) = 0$, which is the case for Fermat type Calabi-Yau hypersurfaces. The difficulty to generalize the construction for these two types of regions to the rest of $\pi^{-1}(\Delta_0^{(n-1)})$ in general is the need to deal with ``thin torus" fibres that have small circles in some directions and big circles in other directions.\\

In this paper, we will construct generalized special Lagrangian fibration for $X_t$ over two other types of regions in $\partial \Delta$. The first type is in the neighborhood of $\pi^{-1}(\Delta_0^{(0)})$. The second type is in the neighborhood of $\pi^{-1}(\Delta_0^{(n-1)}\setminus \Delta_0^{(0)})$. Our construction for the thin torus in regions of the first type is rather complete. Notice that the mirror of quintic Calabi-Yau satisfies $\Delta_0^{(n-1)}=\Delta_0^{(0)}$. Therefore, this construction will be enough to produce monodromy representing generalized special Lagrangian fibration for the mirror family of quintic Calabi-Yau. Our construction for the thin torus in regions of the second type is less satisfactory, since it will require stronger sense of near the large complex limit for the hypersurfaces and toric metric. Nerverthless, this construction will produce monodromy representing generalized special Lagrangian fibration for such Calabi-Yau hypersurfaces with suitable toric metric.\\

In the sequels of this paper, we will give more satisfactory treatment of the second type thin torus and discuss some singular fibres.\\

This paper is organized as follows. In section 2 we lay out the local setting that each of our constructions in section 5 and 6 will reduce to and clarify the estimates for the thin torus needed to apply the implicit function theorem. These estimates are worked out in section 3. In section 4, applying the estimates of thin torus from section 3, we work out a precise quantitative version of theorem 5.1 in \cite{sl1}, which will be needed for the construction in section 6. Sections 5 and 6 are devoted to the construction of generalized special Lagrangian fibrations over the two types of regions based on the estimates from section 3. As consequence, we are able to construct monodromy representing generalized special Lagrangian torus fibration for a large class of Calabi-Yau hypersurfaces in toric variety including the quintic Calabi-Yau 3-fold (corollary \ref{fc}) and its mirror (corollary \ref{ec}). We also indicated how to generalize our construction to general Calabi-Yau hypersurfaces in toric variety in a remark in section 6.\\

\se{Basic setting}
Using suitable toric coordinate, locally around the region we are interested, we can express the hypersurface $X_t$ as a hypersurace in $\mathbb{C}^{n+1}$:

\[
X_t = \{z\in \mathbb{C}^{n+1}|\tilde{p}_t(z)=0\},\ \ {\rm where}\ \tilde{p}_t(z) = z^m + tp(z),\ z^m = \prod_{k=0}^n z_k^{m_k},
\]

$\{0,\cdots,n\} = I' \cup I''$, $m_k\leq 0$ for $k\in I'$ and $m_k=1$ for $k\in I''$. Correspondingly we have decompositions $m = (m_0,\cdots, m_n) = (m',m'')$ and $z = (z_0,\cdots, z_n) = (z',z'')$. We will start with the local model.\\

{\bf Local model:} In $\mathbb{C}^{n+1}$, consider the family of hypersurfaces $\{X_{t,0}\}$ defined as

\[
X_{t,0} = \{z^m =t\}.
\]

Clearly $X_{t,0}$ are all toric varieties. For a toric metric on $\mathbb{C}^{n+1}$ with \k form $\omega$, $\omega_{t,0} = \omega|_{X_{t,0}}$ is a toric \k form on $X_{t,0}$. Therefore the family of natural real torus
 
\[
L_{t,0} = \{z\in X_{t,0}||z_k|=r_k ({\rm constant}),\ {\rm for}\ 0\leq k \leq n\}
\]

form the generalized special Lagrangian torus fibration of $X_{t,0}$ with respect to $\omega_{t,0}$. The fibration can also be defined using the fibration map $\hat{\pi}_{t,0}: X_{t,0} \rightarrow \partial \Delta$ defined as $\displaystyle \hat{\pi}_{t,0}(z) = \{|z_k|^2 - r_{\min}^2\}_{k=0}^n$, where $\Delta$ denotes the first quadrant in $\mathbb{R}^{n+1}$ and $\displaystyle r_{\min} = \min_{0\leq k\leq n} |z_k|$.
\begin{flushright} \rule{2.1mm}{2.1mm} \end{flushright}
Let $\check{p}(z) = \log p(z)$, and

\[
X_{t,s} = \{\tilde{p}_{t,s}(z)=0\},\ \ {\rm where}\ \tilde{p}_{t,s}(z) = z^m + te^{s\check{p}(z)}.
\]

We will use the family $\{X_{t,s}\}_{s\in [0,1]}$ to connect the local model $X_{t,0}$ and $X_t = X_{t,1}$. Our idea here is to deform the generalized special Lagrangian torus fibration of the local model $X_{t,0}$ to the generalized special Lagrangian torus fibration of $X_t$.\\

Hamiltonian-gradient vector fields were discussed in section 4 of \cite{sl1}. Here we will need explicit expressions of them. Generally, consider complex hypersurfaces $Y_t$ defined by $t = u/v$, where $u$ and $v$ are holomorphic functions. $t$ is meromorphic when the family $\{Y_t\}$ has non-empty base locus. Let $V$ be the Hamiltonian-gradient vector field with respect to the family $\{Y_t\}$.\\
\begin{lm}
\label{bc}
When restricted to the complex hypersurface $Y_t$, where $t$ is a constant, we have

\[
V = \frac{\nabla f}{|\nabla f|} = \frac{{\rm Re}(\bar{v}\nabla (u-tv))}{|d(u-tv)|^2},\ \ {\rm where}\ f = {\rm Re}(t).
\]

In particular, when $Y_t$ are all smooth, $V$ will be smooth vector field on the total space ${\cal Y}$.
\end{lm}

{\bf Proof:} For $t=u/v$, we have

\[
dt = \frac{du-tdv}{v}.
\]

When restricted to the complex hypersurface $Y_t$, where $t$ is a constant, we have

\[
\nabla t = \frac{1}{|v|^2}\bar{v}\nabla(u-tv).
\]

\[
\nabla f = \frac{1}{|v|^2}{\rm Re}(\bar{v}\nabla(u-tv)).
\]

\[
|\nabla f|^2 = |df|^2 = \frac{1}{2}|dt|^2 = \frac{1}{2}\frac{|du-tdv|^2}{|v|^2}.
\]

\[
V = \frac{\nabla f}{|\nabla f|} = \frac{2{\rm Re}(\bar{v}\nabla (u-tv))}{|d(u-tv)|^2}.
\]
\begin{flushright} \rule{2.1mm}{2.1mm} \end{flushright}
In lemma \ref{bc}, we can take $u = \log(-z^m/t)$, $v = \check{p}$ and replace $t$ by $s$, $Y_t$ by $X_{t,s}$. Then along $X_{t,s}$,

\[
d(u-sv) = \sum_{k=0}^n q_k\frac{dz_k}{z_k},\ \ {\rm where}\ q_k = m_k - sz_k\check{p}_k.
\]

The normalized Hamiltonian-gradient vector field for the parameter $s$ is

\[
V = 2{\rm Re}\left(\check{p} \sum_{i=0}^n \lambda_i z_i\frac{\partial}{\partial z_i}\right),\ \ {\rm where}
\]
\begin{equation}
\label{lambda}
\lambda_i = \left(\sum_{j=0}^n\frac{g^{i\bar{j}}\bar{q}_j}{z_i\bar{z}_j}\right) \left(\sum_{k,j=0}^n\frac{g^{k\bar{j}}q_k \bar{q}_j}{z_k\bar{z}_j} \right)^{-1}
\end{equation}

The normalized Hamiltonian-gradient vector field for the parameter $t$ is

\[
V_t = 2{\rm Re}\left(\frac{1}{t} \sum_{i=0}^n \lambda_i z_i\frac{\partial}{\partial z_i}\right).
\]

Let

\[
\Omega = \frac{1}{p}\bigwedge_{i=0}^n dz_i,\ \ \Omega_t = \imath(V_t)\Omega|_{X_{t,s}},\ \ {\rm where}\ p=e^{s\check{p}}.
\]

(Here $\Omega$, $\Omega_t$, $V_t$ and $V$ are also depending on $s$. We are omitting the subscript for the simplicity of notations.) Notice that according to lemma 4.1 of \cite{sl1}, $\Omega_t$ so defined when $s=1$ coincides with the holomorphic volume form $\Omega_t$ on $X_t = X_{t,1}$. Let $\phi_s$ denote the flow for $V$. For fixed $t$, the family $(X_{t,s},\omega_{t,s},\Omega_t)$ of Calabi-Yau manifolds can be reduced to the equivalent family of Calabi-Yau structures $\{\Omega_{t,s} = \phi_s^*\Omega_t\}$ on the fixed symplectic manifold $(X_{t,0},\omega_{t,0})$.\\

For the generalized special Lagrangian submanifold $L_{t,0}$ in $(X_{t,0},\omega_{t,0},\Omega_{t,0})$, there is an identification $X_{t,0}\cong T^*L_{t,0}$ near $L_{t,0}$ such that $\omega_{t,0}$ is identified with the canonical symplectic form on $T^*L_0$. (See \cite{W} for the general case. In our toric situation, such identification is actually explicit and does not need the general result from \cite{W}.) Any Lagrangian submanifold $L$ near $L_{t,0}$ can be identified with the graph of $dh$ for some smooth function $h$ on $L_{t,0}$. Under standard coordinate $(x,y)$ on $T^*L_{t,0}$, $\omega_{t,0} = dx\wedge dy$ and $L$ is locally the graph $(x, \frac{\partial h}{\partial x})$.\\

Locally, (as in \cite{sl1}, with slight adjustment of notation,) we can write

\[
\Omega_{t,s} = \eta_s\bigwedge_{k=1}^n (dx_k + u_{s,k}dy),
\]  
\[
\Omega_{t,s}|_L = \eta_s\left(x, \frac{\partial h}{\partial x}\right)\det\left(I + U_s\frac{\partial^2 h}{\partial x^2}\right)\bigwedge_{k=1}^n dx_k,
\] 
where $U_s = (u_{s,1},\cdots,u_{s,n}) = (u_{s,jk})$ and $u_{s,k}dy = \displaystyle\sum_{j=1}^n u_{s,jk}dy_j$.

\[
F(h,s) = {\rm Im}\left(\log \Omega_{t,s}|_L\right) = {\rm Im}\left(\log \eta_s\left(x, \frac{\partial h}{\partial x}\right)+ \log \det\left(I + U_s\frac{\partial^2 h}{\partial x^2}\right)\right)
\] 

defines a map $F: {\cal B}_1\times \mathbb{R} \rightarrow {\cal B}_2$, where ${\cal B}_1 = C^{2,\alpha}(L_{t,0})$ and ${\cal B}_2 = C^\alpha(L_{t,0})$ are Banach spaces. We intend to apply implicit function theorem to $F$ to construct the family of generalized special Lagrangians $L_{t,s}$ with respect to $(X_{t,0},\omega_{t,0},\Omega_{t,s})$.\\

Straightforward computations give us

\[
\delta_h F(h,s) = \frac{\partial F}{\partial h}(h,s) \delta h = a_s^{ij}\frac{\partial^2 \delta h}{\partial x_i\partial x_j} + b_s^i\frac{\partial\delta h}{\partial x_i},
\]

where 

\begin{equation}
\label{bb}
a_s^{ij}(x) = {\rm Im}\left(\left(I + U_s\frac{\partial^2 h}{\partial x^2}\right)^{-1} U_s\right)\left(x,\frac{\partial h}{\partial x}\right),
\end{equation}
\[
b_s^i(x) = {\rm Im}\left(\frac{1}{\eta_s}\frac{\partial \eta_s}{\partial y_i} + {\rm Tr}\left(\left(I + U_s\frac{\partial^2 h}{\partial x^2}\right)^{-1}\frac{\partial U_s}{\partial y_i}\frac{\partial^2 h}{\partial x^2}\right)\right)\left(x,\frac{\partial h}{\partial x}\right).
\]

In particular,

\[
\delta_h F(0,0) = \frac{\partial F}{\partial h}(0,0) \delta h = a^{ij}\frac{\partial^2 \delta h}{\partial x_i\partial x_j} + b^i\frac{\partial\delta h}{\partial x_i},
\]

where 

\[
a^{ij}(x) = {\rm Im}\left((U_0)_{ij}\right)(x,0),\ \ b^i(x) = {\rm Im}\left(\frac{1}{\eta_0}\frac{\partial \eta_0}{\partial y_i}\right)(x,0).
\]

In our toric situation, under the coordinate $x_k = \theta_k$ for $1\leq k \leq n$, $g_{t,0}|_{L_{t,0}}$ is a flat metric on $L_{t,0}$. According to (3.1) in \cite{sl1}, we have $a^{ij}(x) = (g_{t,0}|_{L_{t,0}})^{ij}$ being a constant positive definite matrix and $b^i(x) = 0$. Therefore 

\[
\delta_h F(0,0) = \frac{\partial F}{\partial h}(0,0) \delta h = (g_{t,0}|_{L_{t,0}})^{ij}\frac{\partial^2 \delta h}{\partial x_i\partial x_j}
\]

is just the standard Laplace operator on $L_{t,0}$.\\

\begin{prop}
\label{ba}
Assume that the torus $(L_{t,0}, g_{t,0}|_{L_{t,0}})$ has bounded diameter. Then there exists a constant $C$ (only depending on the dimension) such that

\[
\|\delta h\|_{{\cal B}_1} \leq C \left\|\left(\frac{\partial F}{\partial h}(0,0)\right)\delta h\right\|_{{\cal B}_2}.
\]
\end{prop}
{\bf Proof:} One only need to show that $C$ is independent of the thinness of the torus. It is easy to see that there exists a finite covering map $\pi: \tilde{L} \rightarrow L_{t,0}$ such that $(\tilde{L}, \tilde{g} = \pi^*g_{t,0}|_{L_{t,0}})$ is of bounded geometry. We will use those symbols with ``$\ \tilde{ }\ $" to denote the corresponding pullback objects by $\pi$. By standard Schauder estimate, we have

\[
\|\tilde{\delta h}\|_{\tilde{{\cal B}}_1} \leq C \left\|\left(\frac{\partial \tilde{F}}{\partial \tilde{h}}(0,0)\right)\tilde{\delta h}\right\|_{\tilde{{\cal B}}_2},
\]

where the constant $C$ only depends on the dimension. Since $\tilde{\delta h}$ and $\left(\frac{\partial \tilde{F}}{\partial \tilde{h}}(0,0)\right)\tilde{\delta h}$ on $\tilde{L}$ are invariant under the deck transformations, we have

\[
\|\delta h\|_{{\cal B}_1} = \|\tilde{\delta h}\|_{\tilde{{\cal B}}_1} \leq C \left\|\left(\frac{\partial \tilde{F}}{\partial \tilde{h}}(0,0)\right)\tilde{\delta h}\right\|_{\tilde{{\cal B}}_2} = C \left\|\left(\frac{\partial F}{\partial h}(0,0)\right)\delta h\right\|_{{\cal B}_2}.
\]
\begin{flushright} \rule{2.1mm}{2.1mm} \end{flushright}
According to the implicit function theorem (theorem 3.2 in \cite{sl1}), proposition \ref{ba} has reduced the construction of the family of generalized special Lagrangians $L_{t,s}$ with respect to $(X_{t,0},\omega_{t,0},\Omega_{t,s})$ to the estimates of 

\[
\left\|\frac{\partial F}{\partial h}(h,s)-\frac{\partial F}{\partial h}(0,0)\right\|, {\ \rm and}\ \ \|F(0,s)\|_{{\cal B}_2},
\]

which will be discussed in the following section. In the following sections, we will always assume that $\|h\|_{{\cal B}_1} = \|h\|_{C^{2,\alpha}}$ is bounded (independent of $t$), which will also be part of the condition to apply the implicit function theorem.\\

\se{Basic estimates}
In this section, we will derive the basic estimates along the flow of $V$. The parallel estimates along the flow of $V_t$ will also be mentioned in the remarks following the corresponding results of $V$.\\

Let $Z = {\rm Diag}(z_1,\cdots,z_n)$, $G = G^* = (g_{i\bar{j}})_{n\times n}$. $G$ is called {\bf T-bounded} if $ZGZ^{-1}$ is bounded. Notice that as consequence of T-boundedness, we also have $\bar{Z}^{-1}G\bar{Z}$ and $G$ being bounded.\\

\begin{lm}
\label{cd}
If $G$ is T-bounded and $\det (G) \geq c >0$, then $G^{-1}$ is T-bounded.
\end{lm}
{\bf Proof:} $G$ being bounded and $\det (G) \geq c >0$ imply that $G^{-1}$ is bounded. Since $\det (ZGZ^{-1}) = \det (G)$, we also have $(ZGZ^{-1})^{-1} = ZG^{-1}Z^{-1}$ being bounded. Consequently, $G^{-1}$ is T-bounded.
\begin{flushright} \rule{2.1mm}{2.1mm} \end{flushright}
{\bf Remark:} When apply this lemma, it is often convenient to verify the stronger conditions $G \geq cI$ ($c>0$) or $G^{-1}$ being bounded in place of $\det (G) \geq c >0$.\\
 
{\bf Example:} Assume that $\omega_g$ is a toric \k metric on $\mathbb{C}^{n+1}$ with the toric \k potential $\rho$ (as function of $\{|z_k|^2\}_{k=0}^n$). Then

\[
g_{j\bar{k}} = \rho_k\delta_{jk} + \rho_{jk}\bar{z}_jz_k,\ \ {\rm where}\ \rho_k = \frac{\partial \rho}{\partial |z_k|^2},\ \rho_{jk} = \frac{\partial^2 \rho}{\partial |z_j|^2\partial |z_k|^2}.
\]

Let $G = (g_{j\bar{k}})$. Then $\displaystyle \frac{g_{j\bar{k}}z_j}{z_k} = \rho_k\delta_{jk} + \rho_{jk}|z_j|^2$ (which are in fact smooth functions of $\{|z_k|^2\}_{k=0}^n$) are clearly bounded. Therefore, $G$ is T-bounded.
\begin{flushright} \rule{2.1mm}{2.1mm} \end{flushright}
Recall that $\phi_s$ denotes the flow for $V$. Let $(z_0,\cdots,z_n) = \phi_s(z^0_0,\cdots,z^0_n)$ and

\[
\|\check{p}\| = |\check{p}| + \sum_{k=0}^n\left|z_k \frac{\partial \check{p}}{\partial z_k}\right|.
\]

In the following, we will always assume that $|z^0_0|$ is the smallest among $|z^0_k|$ for $0\leq k \leq n$, and use $(z_1,\cdots,z_n)$ ($(z^0_1,\cdots,z^0_n)$) as the local coordinate of $X_{t,s}$ ($X_{t,0}$).\\

\begin{lm}
\label{cb}
\[
\lambda_k = O\left(\left|\frac{z^0_0}{z^0_k}\right|^2\right),\ \frac{\partial \log z_k}{\partial \log z^0_j} = \delta_{jk} + O\left(\|\check{p}\|\left|\frac{z^0_0}{z^0_k}\right|^2\right),\ \frac{\partial \log z_k}{\partial \log \bar{z}^0_j} = O\left(\|\check{p}\|\left|\frac{z^0_0}{z^0_k}\right|^2\right).
\]

More generally, all the multi-derivatives of $\log z_k - \log z^0_k$ or $\lambda_k$ with respect to $(\log z^0,\log \bar{z}^0)$ are of order $O\left(\|\check{p}\|\left|\frac{z^0_0}{z^0_k}\right|^2\right)$ or $O\left(\left|\frac{z^0_0}{z^0_k}\right|^2\right)$.\\
\end{lm}

{\bf Proof:}
According to (\ref{lambda}), 

\[
\lambda_i = \left(\frac{1}{|z_i|^2}\sum_{j=0}^n\bar{q}_j\frac{g^{i\bar{j}}\bar{z}_i}{\bar{z}_j}\right) \left(\sum_{k,j=0}^n\frac{g^{k\bar{j}}q_k \bar{q}_j}{z_k\bar{z}_j} \right)^{-1}
\]

Since the matrix $(g_{j\bar{k}})$ is T-bounded and $(g^{j\bar{k}})$ is bounded, by lemma \ref{cd}, $(g^{j\bar{k}})$ is also T-bounded. Consequently, $\displaystyle \frac{g^{i\bar{j}}\bar{z}_i}{\bar{z}_j}$ is bounded. Hence $\lambda_k = O\left(\left|\frac{z_0}{z_k}\right|^2\right)$. Since $\lambda_k$ are bounded smooth functions of the bounded terms $\displaystyle \left\{z_k, \bar{z}_k, \frac{|z_0|^2}{|z_k|^2}\right\}_{k=1}^n$. It is straightforward to verify that all of the multi-derivatives of $\lambda_k$ with respect to $(\log z,\log \bar{z})$ are also of order $O\left(\left|\frac{z_0}{z_k}\right|^2\right)$.\\

\[
\frac{d}{ds}\log z_k = V(\log z_k) = \check{p} \lambda_k.
\]

\begin{eqnarray*}
\frac{d}{ds}\frac{\partial \log z_k}{\partial \log z^0_j} &=& \frac{\partial \check{p} \lambda_k}{\partial \log z_i}\frac{\partial \log z_i}{\partial \log z^0_j}+\frac{\partial \check{p} \lambda_k}{\partial \log z_0}\frac{\partial \log z_0}{\partial \log z^0_j} + \frac{\partial \check{p} \lambda_k}{\partial \log \bar{z}_i}\frac{\partial \log \bar{z}_i}{\partial \log z^0_j}+\frac{\partial \check{p} \lambda_k}{\partial \log \bar{z}_0}\frac{\partial \log \bar{z}_0}{\partial \log z^0_j}\\
&=& \left(\frac{\partial \check{p} \lambda_k}{\partial \log z_i}- \frac{q_i}{q_0}\frac{\partial \check{p} \lambda_k}{\partial \log z_0}\right)\frac{\partial \log z_i}{\partial \log z^0_j} + \left(\frac{\partial \check{p} \lambda_k}{\partial \log \bar{z}_i}- \frac{\bar{q}_i}{\bar{q}_0}\frac{\partial \check{p} \lambda_k}{\partial \log \bar{z}_0}\right)\frac{\partial \log \bar{z}_i}{\partial \log z^0_j}.
\end{eqnarray*}

\[
\frac{d}{ds}\frac{\partial \log z_k}{\partial \log \bar{z}^0_j} = \left(\frac{\partial \check{p} \lambda_k}{\partial \log z_i}- \frac{q_i}{q_0}\frac{\partial \check{p} \lambda_k}{\partial \log z_0}\right)\frac{\partial \log z_i}{\partial \log \bar{z}^0_j} + \left(\frac{\partial \check{p} \lambda_k}{\partial \log \bar{z}_i}- \frac{\bar{q}_i}{\bar{q}_0}\frac{\partial \check{p} \lambda_k}{\partial \log \bar{z}_0}\right)\frac{\partial \log \bar{z}_i}{\partial \log \bar{z}^0_j}.
\]

Notice that

\[
\frac{\partial \check{p} \lambda_k}{\partial \log z_i}- \frac{q_i}{q_0}\frac{\partial \check{p} \lambda_k}{\partial \log z_0} = O\left(\|\check{p}\|\left|\lambda_k\right|\right) = O\left(\|\check{p}\|\left|\frac{z_0}{z_k}\right|^2\right).
\]

Similarly

\[
\left(\frac{\partial \check{p} \lambda_k}{\partial \log \bar{z}_i}- \frac{\bar{q}_i}{\bar{q}_0}\frac{\partial \check{p} \lambda_k}{\partial \log \bar{z}_0}\right)= O\left(\|\check{p}\|\left|\frac{z_0}{z_k}\right|^2\right).
\]

By standard estimate of linear ODE, we have

\[
\log z_k = \log z^0_k + O\left(\|\check{p}\|\left|\frac{z_0}{z_k}\right|^2\right),\ \frac{z_k}{z^0_k} = 1 + O\left(\|\check{p}\|\left|\frac{z_0}{z_k}\right|^2\right),
\]
\[
\frac{\partial \log z_k}{\partial \log z^0_j} = \delta_{jk} + O\left(\|\check{p}\|\left|\frac{z_0}{z_k}\right|^2\right),\ \frac{\partial \log z_k}{\partial \log \bar{z}^0_j} = O\left(\|\check{p}\|\left|\frac{z_0}{z_k}\right|^2\right).
\]

Consequently

\[
\frac{z^0_0}{z^0_k} = \frac{z_0}{z_k} + O\left(\|\check{p}\|\left|\frac{z_0}{z_k}\right|\right),\ \ O\left(\left|\frac{z_0}{z_k}\right|^2\right) = O\left(\left|\frac{z^0_0}{z^0_k}\right|^2\right).
\]

Since all the multi-derivatives of $\lambda_k$ with respect to $(\log z,\log \bar{z})$ are of order $O\left(|\frac{z_0}{z_k}|^2\right)$, by induction and similar ODE estimate, all the multi-derivatives of $\log z_k - \log z^0_k$ or $\lambda_k$ with respect to $(\log z^0,\log \bar{z}^0)$ are of order $O\left(\|\check{p}\|\left|\frac{z^0_0}{z^0_k}\right|^2\right)$ or $O\left(\left|\frac{z^0_0}{z^0_k}\right|^2\right)$.
\begin{flushright} \rule{2.1mm}{2.1mm} \end{flushright}
{\bf Remark:}
Along the flow of $V_t$, we have similar estimate as in lemma \ref{cb}. In the proof, one need to replace $\|\check{p}\|$ by $\frac{1}{|t|}$ and parameter $s$ by $t$. When $\left|\frac{z_0}{z_k}\right|^2 = O(t^{\alpha'})$ for a constant $\alpha'>0$, we have the conclusion

\[
\lambda_k = O\left(\left|\frac{z_0}{z_k}\right|^2\right),\ \frac{\partial \log z_k}{\partial \log z^0_j} = \delta_{jk} + O\left(\left|\frac{z_0}{z_k}\right|^2\right),\ \frac{\partial \log z_k}{\partial \log \bar{z}^0_j} = O\left(\left|\frac{z_0}{z_k}\right|^2\right).
\]

More generally, all the multi-derivatives of $\log z_k - \log z^0_k$ or $\lambda_k$ with respect to $(\log z^0,\log \bar{z}^0)$ are of order $O\left(\left|\frac{z_0}{z_k}\right|^2\right)$.\\

{\bf Remark on terminology:} In this section, a bounded term or $O(1)$ quite often means a bounded smooth function of the bounded terms $\displaystyle \left\{z^0_k, \bar{z}^0_k, \frac{|z^0_0|^2}{|z^0_k|^2}\right\}_{k=1}^n$. It is easy to see that the derivatives of such function with respect to $\displaystyle \left\{\log|z^0_k|^2\right\}_{k=1}^n$ will be functions of the same type, therefore are also bounded or $O(1)$.
\begin{flushright} \rule{2.1mm}{2.1mm} \end{flushright}
Let $\displaystyle\omega_{t,0} = \omega|_{X_{t,0}} = \partial \bar{\partial} \tilde{\rho} = \sum_{i,j=1}^n\tilde{g}_{i\bar{j}}dz^0_i d\bar{z}^0_j$ and $\displaystyle\omega^0_{t,0} = \sum_{k=1}^n dz^0_kd\bar{z}^0_k$, where $\tilde{\rho} = \rho|_{X_{t,0}}$. Then\\ 
\begin{lm}
\label{ce}
There exists a constant $c_1>0$ such that 

\[
c_1 \omega^0_{t,0} \leq \omega_{t,0} \leq c_1^{-1} \omega^0_{t,0}.
\]

Both $(\tilde{g}_{i\bar{j}})_{n\times n}$ and its inverse matrix $(\tilde{g}^{i\bar{j}})_{n\times n}$ are T-bounded.
\end{lm}
{\bf Proof:} Clearly $\displaystyle\omega^0 = \sum_{k=0}^n dz^0_kd\bar{z}^0_k$ is quasi-isometric to $\omega$. When restricted to $X_{t,0}$, 

\[
dz^0_0 = -\sum_{k=1}^n m_k\frac{z^0_0}{z^0_k} dz^0_k.
\]

Since we are considering the region where $|z^0_0|$ is the smallest among all $|z^0_k|$, $\omega^0_{t,0}$ is quasi-isometric to $\omega^0|_{X_{t,0}}$, which is quasi-isometric to $\omega_{t,0} = \omega|_{X_{t,0}}$.\\ 
\[
\tilde{g}_{i\bar{j}} = g_{i\bar{j}} - g_{0\bar{j}}\frac{m_iz^0_0}{z^0_i} - g_{i\bar{0}}\frac{m_j\bar{z}^0_0}{\bar{z}^0_j} + g_{0\bar{0}}\frac{m_im_j|z^0_0|^2}{z^0_i\bar{z}^0_j}.
\]

It is straightforward to check that $(\tilde{g}_{i\bar{j}})_{n\times n}$ is T-bounded based on the T-boundedness of $(g_{i\bar{j}})_{(n+1)\times (n+1)}$. Then by lemmas \ref{cd} and the first part of this lemma, $(\tilde{g}^{i\bar{j}})_{n\times n}$ is also T-bounded.
\begin{flushright} \rule{2.1mm}{2.1mm} \end{flushright}
Notice that

\[
\omega_{t,0} = \omega|_{X_{t,0}} = \sum_{k=0}^n d\theta^0_k d(|z^0_k|^2\rho_k) = \sum_{k=1}^n d\theta^0_k d(|z^0_k|^2\rho_k - m_k|z^0_0|^2\rho_0).
\]

When we choose $x_k = \theta^0_k$, we may take $y_k = |z^0_k|^2\tilde{\rho}_k - C_k = |z^0_k|^2\rho_k - m_k|z^0_0|^2\rho_0 - C_k$ for suitable constants $C_k$, $1\leq k \leq n$, so that $\displaystyle\omega_{t,0} = \sum_{k=1}^ndx_kdy_k$ and $y|_{L_{t,0}} =0$. We have

\begin{lm}
\label{cf}
\[
\frac{\partial \log |z^0_j|^2}{\partial y_k} = \min\left(\frac{1}{|z^0_j|^2},\frac{1}{|z^0_k|^2}\right)O(1).
\]

Consequently

\[
\frac{\partial f}{\partial y_k} = \frac{1}{|z^0_k|^2}O\left(\frac{\partial f}{\partial \log |z^0|^2}\right).
\]
\end{lm}
{\bf Proof:}
\[
\frac{\partial y_j}{\partial \log |z^0_k|^2} = \tilde{g}_{j\bar{k}}z^0_j\bar{z}^0_k.
\]

\[
\frac{\partial \log |z^0_j|^2}{\partial y_k} = \frac{1}{z^0_j\bar{z}^0_k}\tilde{g}^{j\bar{k}} = \frac{1}{|z^0_j|^2}\frac{\tilde{g}^{j\bar{k}}\bar{z}^0_j}{\bar{z}^0_k}
= \frac{1}{|z^0_k|^2}\frac{\tilde{g}^{j\bar{k}}z^0_k}{z^0_j}
\]

By lemma \ref{ce}, $(\tilde{g}^{j\bar{k}})_{n\times n}$ is T-bounded. Consequently, $\displaystyle \frac{\tilde{g}^{j\bar{k}}\bar{z}^0_j}{\bar{z}^0_k}$ and $\displaystyle \frac{\tilde{g}^{j\bar{k}}z^0_k}{z^0_j}$ are bounded, which implies the first claim of the lemma. Then

\[
\frac{\partial f}{\partial y_k} = \frac{\partial f}{\partial \log |z^0_j|^2}\frac{\partial \log |z^0_j|^2}{\partial y_k} = \frac{1}{|z^0_k|^2}O\left(\frac{\partial f}{\partial \log |z^0|^2}\right).
\]
\begin{flushright} \rule{2.1mm}{2.1mm} \end{flushright}
For a function $f$ on $\mathbb{R}^n$, we will use $[f]_{C^{k,\alpha}}$ to denote the $(k,\alpha)$-portion of the $C^{k,\alpha}$-norm of $f$. For example

\[
[f]_{C^{1,\alpha}} = \sup_{\tiny{\begin{array}{c}1\leq i \leq n\\x,y \in \mathbb{R}^n\end{array}}} \left(\frac{|\partial_i f(x) - \partial_i f(y)|}{|x-y|^\alpha}\right),\ \ [f]_{C^1} = \sup_{\tiny{\begin{array}{c}1\leq i \leq n\\x \in \mathbb{R}^n\end{array}}} \left(|\partial_i f(x)|\right).
\]
 
Let $h$ be a function on $L_{t,0}$, we have the following estimates.\\

\begin{lm}
\label{cg}
\[
\left[\frac{\partial h}{\partial x_k}\right]_{C^0} \leq C(r_k^0)^{2+\alpha}[h]_{C^{2,\alpha}},\ \left[\frac{\partial^2 h}{\partial x_j\partial x_k}\right]_{C^0} \leq Cr^0_jr^0_k\min((r^0_j)^\alpha,(r^0_k)^\alpha)[h]_{C^{2,\alpha}}.
\]
\[
\left[\frac{\partial h}{\partial x_k}\right]_{C^\alpha} \leq C(r^0_k)^2[h]_{C^{2,\alpha}},\ \left[\frac{\partial^2 h}{\partial x_j\partial x_k}\right]_{C^\alpha} \leq Cr^0_jr^0_k[h]_{C^{2,\alpha}}.
\]
\end{lm}
{\bf Proof:} We will prove the most challenging estimate

\[
\left[\frac{\partial h}{\partial x_k}\right]_{C^\alpha} \leq C(r^0_k)^2[h]_{C^{2,\alpha}}.
\]

When $|x-x'|\geq 10r^0_k$, we have

\[
\frac{1}{|x-x'|^\alpha}\left|\frac{\partial h}{\partial x_k}(x) - \frac{\partial h}{\partial x_k}(x')\right| \leq \frac{1}{|x-x'|^\alpha}\left(\left|\frac{\partial h}{\partial x_k}(x)\right| + \left|\frac{\partial h}{\partial x_k}(x')\right|\right) 
\]
\[
\leq \frac{C}{(r^0_k)^\alpha}\left[\frac{\partial h}{\partial x_k}\right]_{C^0} \leq \frac{C}{(r^0_k)^\alpha}(r_k^0)^{2+\alpha}[h]_{C^{2,\alpha}} \leq C(r_k^0)^2[h]_{C^{2,\alpha}}.
\]

Rescaling the coordinate and metric as $\hat{x} = x/r^0_k$, $\hat{g} = g/(r^0_k)^2$. Then the estimate is reduced to

\[
\left[\frac{\partial h}{\partial \hat{x}_k}\right]_{C^\alpha} \leq C[h]_{C^{2,\alpha}}.
\]

It is easy to check that this estimate is always true in a region with bounded diameter. In our remaining case, $|x-x'|\leq 10r^0_k$ or equivalently $|\hat{x}-\hat{x}'|\leq 10$ ensures that the region of discussion is bounded with respect to the rescaled metric. The other estimates are similar and easier to prove.
\begin{flushright} \rule{2.1mm}{2.1mm} \end{flushright}
\[
\frac{dz_i}{z_i} = \frac{dr_i}{r_i} + id\theta_i = \left(\frac{\partial \log r_i}{\partial \log r^0_k} + i\frac{\partial \theta_i}{\partial \log r^0_k}\right)\frac{dr^0_k}{r^0_k} + \left(\frac{\partial \log r_i}{\partial \theta^0_k} + i\frac{\partial \theta_i}{\partial \theta^0_k}\right)d\theta^0_k.
\]

\[
\Omega_t = \imath(V_t)\Omega|_{X_t} = -\frac{1}{q_0}\prod_{i=1}^n \frac{dz_i}{z_i}\]
\[
= \eta_s\prod_{i=1}^n \left(d\theta^0_i + \left(\frac{\partial \log r_j}{\partial \theta^0_i} + i\frac{\partial \theta_j}{\partial \theta^0_i}\right)^{-1}\left(\frac{\partial \log r_j}{\partial \log r^0_k} + i\frac{\partial \theta_j}{\partial \log r^0_k}\right)\frac{dr^0_k}{r^0_k}\right).
\]

\[
\eta_s = -\frac{1}{q_0}\det\left(\frac{\partial \log r_i}{\partial \theta^0_k} + i\frac{\partial \theta_i}{\partial \theta^0_k}\right).
\]
\[
(U_s)_{ij} = \left(\frac{\partial \log r_l}{\partial \theta^0_i} + i\frac{\partial \theta_l}{\partial \theta^0_i}\right)^{-1}\left(\frac{\partial \log r_l}{\partial \log r^0_k} + i\frac{\partial \theta_l}{\partial \log r^0_k}\right) \frac{\partial \log r^0_k}{\partial y_j}. 
\]

\begin{lm}
\label{ca}
\[
\left|\log |z^0_k|^2 \left(x,\frac{\partial h}{\partial x}\right) - \log |z^0_k|^2 \left(x,0\right)\right| \leq C[h]_{C^2}.
\] 
\end{lm}
{\bf Proof:} With the help of lemma \ref{cf} and \ref{cg}, we have
\[
\left|\log |z^0_k|^2 \left(x,\frac{\partial h}{\partial x}\right) - \log |z^0_k|^2 \left(x,0\right)\right| \leq \sum_{j=1}^n \left|\frac{\partial \log |z^0_k|^2}{\partial y_j}\right|\left|\frac{\partial h}{\partial x_j}\right|
\] 
\[
\leq C\max_{1\leq j\leq n} \left(\frac{1}{|z^0_j|^2}\left|\frac{\partial h}{\partial x_j}\right|\right) \leq C[h]_{C^2}.
\]
\begin{flushright} \rule{2.1mm}{2.1mm} \end{flushright}
{\bf Remark:} According to lemma \ref{ca}, $\displaystyle O\left(|z^0_k|^2 \left(x,\frac{\partial h}{\partial x}\right)\right) = O\left(|z^0_k|^2 \left(x,0\right)\right)$ when $[h]_{C^2}$ is bounded. For this reason, in the following arguments, we will often omit the mentioning of such implicit dependence of $|z^0_k|^2$ on $h$.\\

Let $\displaystyle \xi_\beta = \max_{1\leq l\leq n} \left(\frac{|z^0_0|^2}{|z^0_l|^{2+\beta}}\right)$ for $\beta = 0$ or $\alpha$ (for the flow of $V_t$, replace $z^0$ by $z$ in $\xi_\beta$). We have\\

\begin{lm}
\label{ch}
For $\beta = 0$ or $\alpha$,
\[
\left[\eta_s \left(x,\frac{\partial h}{\partial x}\right) - \eta_0 \left(x,0\right)\right]_{C^{\beta}} \leq C(\|\check{p}\|\xi_\beta + (r_0^0)^{1-\beta}).
\] 
\[
\left[\frac{\partial \eta_s}{\partial y_k} \left(x,\frac{\partial h}{\partial x}\right)\right]_{C^{\beta}} \leq \frac{C}{(r^0_k)^2}(\|\check{p}\|\xi_\beta + (r_0^0)^{1-\beta}).
\]
\[
\left[(U_s)_{ij} \left(x,\frac{\partial h}{\partial x}\right) - (U_0)_{ij} \left(x,0\right)\right]_{C^{\beta}} \leq \frac{C}{r^0_ir^0_j}\left(\|\check{p}\|\xi_\beta + [h]_{C^{2,\beta}}\right).
\] 
\[
\left[(\frac{\partial U_s}{\partial y_k})_{ij} \left(x,\frac{\partial h}{\partial x}\right) - (\frac{\partial U_0}{\partial y_k})_{ij} \left(x,0\right)\right]_{C^{\beta}} \leq \frac{C}{(r^0_k)^2r^0_ir^0_j}\left(\|\check{p}\|\xi_\beta + [h]_{C^{2,\beta}}\right).
\]
(When $\beta=0$, $C^{2,\beta}$ will mean $C^2$. Notice that the second term on the left side of each of the estimates is constant.) 
\end{lm}
{\bf Proof:} The proofs for the four estimates are straightforward and somewhat tedious applications of lemmas \ref{cb} and \ref{cf}. We will carefully prove the third estimate here and the proof of the other three are very similar.\\

We may write $U_s = U^1_s U^2_s$, where

\[
(U^1_s)_{ik} = \left(\frac{\partial \log r_l}{\partial \theta^0_i} + i\frac{\partial \theta_l}{\partial \theta^0_i}\right)^{-1}\left(\frac{\partial \log r_l}{\partial \log r^0_k} + i\frac{\partial \theta_l}{\partial \log r^0_k}\right),\ (U^2_s)_{kj} = \frac{\partial \log r^0_k}{\partial y_j}. 
\]

Estimates in lemma \ref{cb} imply that

\[
\left[(U^1_s)_{ij} \left(x,\frac{\partial h}{\partial x}\right) - (U^1_0)_{ij} \left(x,0\right)\right]_{C^0} \leq C\|\check{p}\|\frac{(r_0^0)^2}{(r_i^0)^2} \leq C\|\check{p}\|\xi_0.
\] 

\[
\left[(\frac{\partial U^1_s}{\partial x})_{ij} \left(x,\frac{\partial h}{\partial x}\right)\right]_{C^0} \leq C\|\check{p}\|\frac{(r_0^0)^2}{(r_i^0)^2} \leq C\|\check{p}\|\xi_0.
\]

\[
\left[(\frac{\partial U^1_s}{\partial y_k})_{ij} \left(x,\frac{\partial h}{\partial x}\right)\right]_{C^0} \leq \frac{C}{(r^0_k)^2}\|\check{p}\|\frac{(r_0^0)^2}{(r_i^0)^2} \leq \frac{C}{(r^0_k)^2}\|\check{p}\|\xi_0.
\]

(The last estimate also need lemma \ref{cf}.) For $\alpha>0$, we have

\[
\left[(U^1_s)_{ij} \left(x,\frac{\partial h}{\partial x}\right) - (U^1_0)_{ij} \left(x,0\right)\right]_{C^\alpha} = \left[(U^1_s)_{ij} \left(x,\frac{\partial h}{\partial x}\right) \right]_{C^\alpha}
\]
\[
\leq C\left[(\frac{\partial U^1_s}{\partial x})_{ij} \left(x,\frac{\partial h}{\partial x}\right)\right]_{C^0} \max_{1\leq l\leq n} \left(|z^0_l|^{-\alpha}\right)
+ C\max_{1\leq k\leq n} \left(\left[(\frac{\partial U^1_s}{\partial y_k})_{ij} \left(x,\frac{\partial h}{\partial x}\right)\right]_{C^0}\left[\frac{\partial h}{\partial x_k}\right]_{C^{\alpha}}\right) 
\]
\[
\leq C\left(\|\check{p}\|\frac{(r_0^0)^2}{(r_i^0)^2}\max_{1\leq l\leq n} \left(|z^0_l|^{-\alpha}\right) + \frac{(r_0^0)^2}{(r_i^0)^2}[h]_{C^{2,\alpha}}\right)\leq C\left(\|\check{p}\|\xi_\alpha + [h]_{C^{2,\alpha}}\right).
\] 

Estimates in lemma \ref{cb} imply that

\[
\left[(U^2_s)_{ij} \left(x,\frac{\partial h}{\partial x}\right)\right]_{C^0} 
 \leq \frac{C}{(r^0_j)^2}.
\]

Since $U^2_s$ is independent of $s$ and $x$, we have

\[
\left[(U^2_s)_{ij} \left(x,\frac{\partial h}{\partial x}\right) - (U^2_0)_{ij} \left(x,0\right)\right]_{C^\beta} 
\]
\[
\leq C\max_{1\leq k\leq n} \left(\left[(\frac{\partial U^2_s}{\partial y_k})_{ij} \left(x,\frac{\partial h}{\partial x}\right)\right]_{C^0}\left[\frac{\partial h}{\partial x_k}\right]_{C^{\beta}}\right) \leq \frac{C}{r^0_ir^0_j}[h]_{C^{2,\beta}}
\]

for $\beta = 0$ or $\alpha$. Combining the estimates for $U^1_s$ and $U^2_s$ gives the estimate for $U_s$. One subtle point is that the estimate of $(U^1_s)_{ik}$ involves the factor $\frac{(r_0^0)^2}{(r_i^0)^2}$ except one term in $[(U^1_s)_{ik}]_{C^0}$ involves $\delta_{ik}$. The estimate of $(U^2_s)_{kj}$ involves the factor $\frac{1}{r^0_kr^0_j}$ (see lemma \ref{cf}). To get the factor $\frac{1}{r^0_ir^0_j}$ in both cases for the final estimate of $(U_s)_{ij}$, we use the following

\[
\frac{(r_0^0)^2}{(r_i^0)^2} \frac{1}{r_k^0r_j^0} = \frac{1}{r^0_ir^0_j}\frac{(r_0^0)^2}{r^0_ir^0_k} \leq \frac{1}{r^0_ir^0_j}\xi_0\ {\rm and}\ \delta_{ik}\frac{1}{r^0_kr^0_j} = \frac{1}{r^0_ir^0_j}.
\]

For the estimates of $\eta_s$, other than those partial derivative terms that can be similarly delt with as in the third estimate, there is also an additional factor $q_0 = 1- sz_0\check{p}_0$ that need the following estimate

\[
[q_0 - 1]_{C^\beta} \leq C(r_0^0)^{1-\beta},
\]

which is straightforward to verify.
\begin{flushright} \rule{2.1mm}{2.1mm} \end{flushright}
Recall the expressions of $a_s^{ij}(x)$ and $b_s^i(x)$ from (\ref{bb}). Applying lemma \ref{ch}, it is easy to derive the following\\

\begin{lm}
\label{ci}
For $\beta = 0$ or $\alpha$,
\[
[a_s^{ij} - a^{ij}]_{C^{\beta}} \leq \frac{C}{r^0_ir^0_j}\left(\|\check{p}\|\xi_\beta + [h]_{C^{2,\beta}}\right).
\]
\[
[b_s^i - b^i]_{C^{\beta}} \leq \frac{C}{(r^0_i)^2}\left(\|\check{p}\|\xi_\beta + [h]_{C^{2,\beta}} + (r_0^0)^{1-\beta}\right).
\]
\end{lm}
{\bf Proof:} Let $\Lambda = {\rm Diag}(r^0_1,\cdots,r^0_n)$, $(a_s^{ij})^\Lambda = \Lambda(a_s^{ij})\Lambda$, $U_s^\Lambda = \Lambda U_s\Lambda$, $\left(\frac{\partial^2 h}{\partial x^2}\right)^\Lambda = \Lambda^{-1}\frac{\partial^2 h}{\partial x^2}\Lambda^{-1}$. Then

\[
(a_s^{ij})^\Lambda = {\rm Im}\left(\left(I + U_s^\Lambda\left(\frac{\partial^2 h}{\partial x^2}\right)^\Lambda\right)^{-1} U_s^\Lambda\right)\left(x,\frac{\partial h}{\partial x}\right),
\]
\[
b_s^i(x) = {\rm Im}\left(\frac{1}{\eta_s}\frac{\partial \eta_s}{\partial y_i} + {\rm Tr}\left(\left(I + U_s^\Lambda\left(\frac{\partial^2 h}{\partial x^2}\right)^\Lambda\right)^{-1}\frac{\partial U_s^\Lambda}{\partial y_i}\left(\frac{\partial^2 h}{\partial x^2}\right)^\Lambda\right)\right)\left(x,\frac{\partial h}{\partial x}\right).
\]

Estimates in lemma \ref{ch} imply that

\[
\left[U_s^\Lambda \left(x,\frac{\partial h}{\partial x}\right) - U_0^\Lambda \left(x,0\right)\right]_{C^{\beta}} \leq C\left(\|\check{p}\|\xi_\beta + [h]_{C^{2,\beta}}\right).
\]

Also notice that $U_0^\Lambda \left(x,0\right)$ is a bounded constant matrix and

\[
\left[\left(\frac{\partial^2 h}{\partial x^2}\right)^\Lambda\right]_{C^{\beta}} = [h]_{C^{2,\beta}}.
\]

Consequently

\[
\left[(a_s^{ij})^\Lambda - (a^{ij})^\Lambda\right]_{C^{\beta}} \leq C\left(\|\check{p}\|\xi_\beta + [h]_{C^{2,\beta}}\right). 
\]

The estimate for $b_s^i$ is similar by applying lemma \ref{ch}.
\begin{flushright} \rule{2.1mm}{2.1mm} \end{flushright}
{\bf Remark:} The corresponding estimates along the flow of $V_t$ are

\[
[a_t^{ij} - a^{ij}]_{C^{\beta}} \leq \frac{C}{r_ir_j}\left(\xi_\beta + [h]_{C^{2,\beta}}\right).
\]
\[
[b_t^i - b^i]_{C^{\beta}} \leq \frac{C}{(r_i)^2}\left(\xi_\beta + [h]_{C^{2,\beta}} + |z_0|^{1-\beta}\right).
\]
\begin{flushright} \rule{2.1mm}{2.1mm} \end{flushright}
\begin{prop}
\label{cc}
\[
\left\|\frac{\partial F}{\partial h}(h,s)-\frac{\partial F}{\partial h}(0,0)\right\| \leq C\left(\|\check{p}\|\xi_\alpha + [h]_{C^{2,\alpha}} + (r_0^0)^{1-\alpha}\right).
\]
\[
\|F(0,s)\|_{{\cal B}_2} \leq C(\|\check{p}\|\xi_\alpha + (r_0^0)^{1-\alpha}).
\]
\end{prop}
{\bf Proof:} Recall that
\[
\left(\frac{\partial F}{\partial h}(h,s)-\frac{\partial F}{\partial h}(0,0)\right) \delta h = (a_s^{ij} - a^{ij})\frac{\partial^2 \delta h}{\partial x_i\partial x_j} + (b_s^i - b^i)\frac{\partial\delta h}{\partial x_i}
\]

and
\[
F(0,s) = {\rm Im}(\log \eta_s(x,0)).
\]

Lemmas \ref{ci}, \ref{cg} will imply the first estimate. The second estimate is implied by lemma \ref{ch}.
\begin{flushright} \rule{2.1mm}{2.1mm} \end{flushright}
{\bf Remark:} We again have the corresponding estimate along the flow of $V_t$.

\[
\left\|\frac{\partial F}{\partial h}(h,t)-\frac{\partial F}{\partial h}(0,0)\right\| \leq C\left(\xi_\alpha + [h]_{C^{2,\alpha}} + |z_0|^{1-\alpha}\right).
\]
\[
\|F(0,t)\|_{{\cal B}_2} \leq C(\xi_\alpha + |z_0|^{1-\alpha}).
\]
\begin{flushright} \rule{2.1mm}{2.1mm} \end{flushright}

\se{Fibration over top dimensional faces in $\partial \Delta$}
In \cite{sl1} section 5, we constructed the generalized special Lagrangian fibration over $U^{\rm top}_t$, which is in the interior of the top dimensional faces of $\partial \Delta$ (also see \cite{G}). In this section, using our estimates in section 3 for the thin torus, we will give a more quantitative characterization of $U^{\rm top}_t$, which is crucial for our argument in section 6.\\

Locally in the face $r^0_0=0$ of $\partial \Delta$, we may use $(r^0_1,\cdots,r^0_n)$ as local coordinate.

Take $X_0 \subset P_\Delta$ to be the large complex limit, with the natural generalized special Lagrangian fibration $X_0 \rightarrow \partial \Delta$. Smooth part of $X_0$ is a union of top dimensional complex torus fibred over top dimensional faces in $\partial \Delta$.\\ 

\begin{theorem}
\label{da}
There exists a constant $\beta>0$ such that when $\displaystyle\min_{1\leq k \leq n}((r^0_k)^2) \geq Ct^\beta$, the corresponding $L_0$ in theorem 5.1 of \cite{sl1} will be able to deform to $L_t \subset X_t$ that form the generalized special Lagrangian fibration over the open set $U^{\rm top}_t$ in the top dimensional faces of $\partial \Delta$.
\end{theorem}

{\bf Proof:} According to the remark after lemma \ref{cb}

\[
\min_{1\leq k \leq n}((r_k)^2) \geq C\min_{1\leq k \leq n}((r^0_k)^2) \geq Ct^\beta.
\]

When $\beta>0$ is chosen suitably (for example $\beta < \frac{2}{n+2}$), there exists $\alpha'>0$ such that $(\xi_\alpha + |z_0|^{1-\alpha}) \leq Ct^{\alpha'}$. Start with $\|h\|_{{\cal B}_1} \leq C_1$, where $C_1$ is small depending on the constant in proposition \ref{ba}. The estimates in proposition \ref{ba} and the remark after proposition \ref{cc} will enable us to apply the implicit function theorem (theorem 3.2 in \cite{sl1}) to prove the theorem. In particular, $h$ for the actual solution will satisfy $\|h\|_{{\cal B}_1} \leq C_2t^{\alpha'}$.\\

For fixed small $t$, to show $L_t$ form a generalized special Lagrangian fibration when the fibre $L_0$ vary, it is sufficient to show that non-trivial deformation 1-forms of $L_t$ have no zeroes on $L_t$.  It is straightforward to check that the deformation 1-forms of $L_t$ are spanned by $\{dx_k - df_k\}_{k=1}^n$, where $\{f_k\}_{k=1}^n$ are functions on $L_t$ satisfying

\[
a_t^{ij}\frac{\partial^2 f_k}{\partial x_i\partial x_j} + b_t^i\frac{\partial f_k}{\partial x_i} = b_t^k,
\]

which can be rewritten as

\[
a^{ij}\frac{\partial^2 f_k}{\partial x_i\partial x_j} = (a^{ij} - a^{ij}_t)\frac{\partial^2 f_k}{\partial x_i\partial x_j} + (b^i - b^i_t)\frac{\partial f_k}{\partial x_i} + (b^k_t - b^k).
\]

Propositions \ref{ba}, lemma \ref{cg} and the remark after lemma \ref{ci} imply that

\[
|f_k|_{C^{2,\alpha}} \leq C\left(\xi_\alpha + [h]_{C^{2,\alpha}} + |z_0|^{1-\alpha}\right)(|f_k|_{C^{2,\alpha}} + \frac{1}{r_k^2}).
\]

Since $[h]_{C^{2,\alpha}} \leq Ct^{\alpha'}$, we have

\[
|f_k|_{C^{2,\alpha}} \leq \frac{C}{r_k^2}\left(\xi_\alpha + [h]_{C^{2,\alpha}} + |z_0|^{1-\alpha}\right)\leq \frac{Ct^{\alpha'}}{r_k^2}.
\]

Apply lemma \ref{cg}, we have

\[
\left[\frac{\partial f_k}{\partial x_k} \right]_{C^0} \leq Cr_k^{2 + \alpha}|f_k|_{C^{2,\alpha}} \leq Ct^{\alpha'}r_k^{\alpha}.
\]

Hence the $dx_k$-coefficient of $dx_k - df_k$ is of order $1+O(t^{\alpha'})$ which is non-zero when $t$ is small. Consequently $dx_k - df_k$ does not vanish anywhere.
\begin{flushright} \rule{2.1mm}{2.1mm} \end{flushright}

\se{Fibration near $\pi^{-1}(\Delta_0^{(0)})$}
Consider the local situation

\[
X_{t,s} = \{z\in \mathbb{C}^{n+1}|\tilde{p}_{t,s}(z)=0\},\ \ {\rm where}\ \tilde{p}_{t,s}(z) = \prod_{k=0}^n z_k + tp(z),\ \ p(z) = e^{s\check{p}(z)}.
\]

Let $z = (z_0,\cdots, z_n) = (z',z'')$. Correspondingly we have $\{0,\cdots,n\} = I' \cup I''$. Assume that $\check{p}(z',0)=0$. Recall the notations from section 2, we have\\

\begin{prop}
\label{ea}
When $t$ is small enough, for the generalized special Lagrangian thin torus $L_{t,0}$ satisfying 
\[
\max_{k\in I''}r^0_k \leq C\min_{1\leq k \leq n}((r^0_k)^{\alpha}), 
\]
there exists the deformation family of generalized special Lagrangian thin torus $\{L_{t,s}\}_{s\in [0,1]}$ with respect to $(X_{t,0}, \omega_{t,0}, \Omega_{t,s})$. Further more, $L_{t,s}$ is part of a generalized special Lagrangian fibration in $(X_{t,0}, \omega_{t,0}, \Omega_{t,s})$.
\end{prop}
{\bf Proof:} Since $\check{p}(z',0)=0$,

\[
\|\check{p}\| = O\left(\max_{k\in I''}r^0_k \right).
\]

Take $\alpha >0$ smaller if necessary, there exists a constant $\alpha'>0$ such that 

\[
\|\check{p}\|\xi_\alpha + (r_0^0)^{1-\alpha} \leq Ct^{\alpha'}.
\]

Start with $\|h\|_{{\cal B}_1} \leq C_1$, where $C_1$ is small depending on the constant in proposition \ref{ba}. The estimates in proposition \ref{ba} and \ref{cc} will enable us to apply the implicit function theorem (theorem 3.2 in \cite{sl1}) to prove the claim. In particular, $h$ for the actual solution will satisfy $\|h\|_{{\cal B}_1} \leq C_2t^{\alpha'}$. Using similar argument as in theorem \ref{da}, one can show that the non-trivial deformation 1-forms of $L_{t,s}$ have no zeroes, therefore $L_{t,s}$ is part of a generalized special Lagrangian fibration in $(X_{t,0}, \omega_{t,0}, \Omega_{t,s})$.
\begin{flushright} \rule{2.1mm}{2.1mm} \end{flushright}
\begin{theorem}
\label{eb}
When $t$ is small enough, there exists an open neighborhood $U^{(0)}_t$ of $\pi^{-1}(\Delta_0^{(0)}) \subset \partial \Delta$ such that for any generalized special Lagrangian fibre $L_{t,0}$ in $X_{t,0}$ over $U^{(0)}_t$, there exist a smooth family $\{L_{t,s}\}_{s\in [0,1]}$, where $L_{t,s}$ is a generalized special Lagrangian torus in $(X_{t,0}, \omega_{t,0}, \Omega_{t,s})$ that is Hamiltonian equivalent to $L_{t,0}$. When the fibre $L_{t,0}$ varies over $U^{(0)}_t$, $\phi_1(L_{t,1}) \subset X_{t,1} = X_t$ will form a generalized special Lagrangian fibration over $U^{(0)}_t$. This fibration will coincide with the fibration over $U^{\rm top}_t$ (in theorems 5.1 and 5.2 of \cite{sl1} and theorem \ref{da}) on the overlaps.
\end{theorem}
{\bf Proof:}
Generically, a connected component of $\pi^{-1}(\Delta_0^{(0)})$ is either a point or a contractable union of top dimensional faces of $\partial \Delta$. The first case was considered in \cite{sl1}. We will deal with the second case here. Near a vertex of $\pi^{-1}(\Delta_0^{(0)})$, $\pi^{-1}(\Delta_0^{(0)})$ can locally be identified with a union of several codimension 1 faces of the first quadrant in $\mathbb{R}^{n+1}$. Assume that $\sigma_0 = \{(r_0^0)^2 = |z_0^0|^2=0\}$ is one of these faces. Then all the subfaces of $\sigma_0$ are $\sigma_I = \{(r_i^0)^2 = 0|i\in I\}$ for $I\ni 0$. The fibration over a tubular neighborhood of $\sigma_I$ can be discussed using proposition \ref{ea} with $I'' = I$.\\

Without lost of generality, we assume that $r^0_k$ is in ascending order with respect to $k$. When $|I|=2$, $I'' = I = \{0,1\}$, we always have 

\[
\min_{1\leq k \leq n} r^0_k = r^0_1 = \max_{k\in I''}r^0_k,
\]

and proposition \ref{ea} will always apply. Therefore the fibration over interior of $\sigma_0$ will extend over each codimension 1 edge of $\sigma_0$ to adjacent codimension 1 faces of $\partial \Delta$.\\

In general, let

\[
X_{t,0}^I = \{z^0\in X_{t,0}||z^0_k|^2 = \mu \leq |z^0_j|^2,\ \ {\rm for\ any}\ k\in I,\ j\not\in I.\}
\]

\[
U_{t,0}^I = \{z^0\in X_{t,0}|\min_{1\leq k \leq n}((r^0_k)^{\alpha}) \geq \max_{k\in I}r^0_k \leq |z^0_j|^2,\ \ {\rm for\ any}\ j\not\in I.\}
\]

Notice that $X_{t,0}^I \subset U_{t,0}^I$ and proposition \ref{ea} applies to $U_{t,0}^I$. $X_{t,0}^I$ fibres over $\sigma_I$, and $U_{t,0}^I$ fibres over a neighborhood $U_{\sigma_I}$ of $\sigma_I$ in $\partial \Delta$. $\displaystyle U_{t,0}^{\sigma_0} = \bigcup_{I\ni 0} U_{t,0}^I$ fibres over $\displaystyle U^{\sigma_0} = \bigcup_{I\ni 0} U_{\sigma_I}$, which is an open neighborhood of $\sigma_0$. We may take $U^{(0)}_{t,0}$ ($U^{(0)}_t$) to be the union of such neighborhoods $U_{t,0}^{\sigma_0}$ ($U^{\sigma_0}$). $U^{(0)}_t$ is an open neighborhood of $\pi^{-1}(\Delta_0^{(0)})$. Apply proposition \ref{ea}, we have an open set $U^{(0)}_{t,1} \subset X_{t,1}$ and fibration $\hat{\pi}_{t,1}: U^{(0)}_{t,1} \rightarrow U^{(0)}_t$.\\

Let $\sigma$ be a top dimension face of $\partial \Delta$ that is adjacent to $\pi^{-1}(\Delta_0^{(0)})$. What we have just proved is enough to make the open neighborhood $U^{(0)}_t$ of $\pi^{-1}(\Delta_0^{(0)})$ and $U^{\rm top}_t \cap \sigma$ touch. We need to show that the two fibrations coincide on the overlap. From the definitions, it is easy to check that for suitable $\{r_k\}_{k=1}^n$ satisfying $\log r_k$ are bounded for $1\leq k \leq n$, in the local models, $L_{t,0}$ ($L_0$) is a fibre over $U^{(0)}_t$ ($U^{\rm top}_t \cap \sigma$), where

\[
L_{t,0} = \{z\in X_{t,0}||z_k|=r_k ({\rm constant}),\ {\rm for}\ 0\leq k \leq n\},
\]
\[
L_0 = \{z_0=0, |z_k|=r_k ({\rm constant}),\ {\rm for}\ 1\leq k \leq n\}.
\]

Notice that $L_{t,0}$ and $L_0$ are torus of bounded geometry. This is exactly the situation as in the proof of theorem 7.1 of \cite{sl1}. The proof of theorem 7.1 of \cite{sl1} implies that after suitable adjustment of $\{r_k\}_{k=1}^n$ in $L_0$, the deformed torus $L_{t,1}$ coincide with $L_1$.\begin{flushright} \rule{2.1mm}{2.1mm} \end{flushright}
As we mentioned in the introduction, the mirror of quintic Calabi-Yau satisfy $\Delta_0^{(n-1)}=\Delta_0^{(0)}$. Apply theorem \ref{eb}, we have\\

\begin{co}
\label{ec}
For the mirror of quintic Calabi-Yau or more generally for the Calabi-Yau family $\{X_t\}$ in a toric variety $P_\Delta$ satisfying $\Delta_0^{(n-1)}=\Delta_0^{(0)}$ with the restriction of any toric \k metric from $P_\Delta$, when $t$ is small enough, there exists a smooth monodromy representing generalized special Lagrangian torus (partial) fibration for $X_t$.
\end{co}
\begin{flushright} \rule{2.1mm}{2.1mm} \end{flushright}

\se{Fibration near $\pi^{-1}(\Delta_0^{(n-1)}\setminus \Delta_0^{(0)})$}
Starting from \cite{sl1}, our discussion so far is valid for any toric \k metric on $P_\Delta$. The construction in this section is of somewhat different nature, where the precise meaning of near the large complex limit discussed in the introduction will be very essential and the toric \k metric on $P_\Delta$ also has to depend on the Calabi-Yau hypersurface $X_t \subset P_\Delta$.\\

The concept of ``near the large complex limit" in this section is slightly different from the one in the introduction. In this section, $X_t$ is said to be near the large complex limit, if $\tau$ and $t = \tau^{-w_{m_o}}$ are small. Namely, $t$ is not allowed to get small independent of $\tau$ once $w$ is determined.\\

Define $\tilde{w}_m = w_m$ for $m\in \Delta_0\setminus \{m_o\}$ and $\tilde{w}_{m_o}=0$. $\{\tilde{w}_m\}_{m\in \Delta_0}$ is clearly still strictly convex on $\Delta_0$ and determines the same simplicial decomposition of $\Delta_0$ as $\{w_m\}_{m\in \Delta_0}$. Let $\tilde{\Delta}_0$ denote those integral $m$ in the real polyhedron spanned by $\Delta_0$. Let $\tilde{\Delta}_0^{(k)}$ denote the $k$-skeleton of $\tilde{\Delta}_0$. ($\Delta_0$ ($\Delta_0^{(k)}$) can be reinterpreted as containing those integral $m$ in $\tilde{\Delta}_0$ ($\tilde{\Delta}_0^{(k)}$) such that $a_m \not=0$.) It is easy to see that one can extend $\{\tilde{w}_m\}_{m\in \Delta_0}$ to $\{\tilde{w}_m (\geq 0)\}_{m\in \tilde{\Delta}_0}$, which is also strictly convex.\\

{\bf Assumption*:} For any $\tilde{m} \in \Delta_0^{(l)}\setminus \Delta_0^{(l-1)}$ ($l\leq n-1$), we assume that there exists a primitive top dimensional simplex $S (\ni \tilde{m})$ in the simplicial decomposition of $\tilde{\Delta}_0$ determined by $\{\tilde{w}_m (\geq 0)\}_{m\in \tilde{\Delta}_0}$ such that $\dim S\cap \tilde{\Delta}_0^{(l)} = l$. We also assume that $P_\Delta$ is equivalent to its anti-canonical model $P_{\Delta_0}$.\\

The following argument will apply to the situation when assumption* is satisfied. It will include most notably the situation of Calabi-Yau quintics in $\mathbb{CP}^4$ with a primitive or a standard simplicial decomposition of $\Delta_0$. (By standard simplicial decomposition, we mean any simplicial subdivision of the standard polyhedron decomposition of $\Delta_0$ resulting from cutting by integral hyperplanes that are parallel to the faces of $\Delta_0$.) We will indicate in a remark at the end of this section what is needed to generalize the discussion to the general case of Calabi-Yau hypersurfaces in toric varieties.\\

By adjustment of an affine function, we can make $\tilde{w}_{m} =0$ for $m\in S$, $\tilde{w}_{m} >0$ for $m\in \tilde{\Delta}_0 \setminus S$. By further adjustment of a suitable affine function, we can make $\tilde{w}_{\tilde{m}}=0$, $\tilde{w}_{m} =\delta$ for $m\in S \setminus \{\tilde{m}\}$, $\tilde{w}_{m} > C\delta$ for $m\in \tilde{\Delta}_0 \setminus S$, where $\delta>0$ is a small constant and $C>0$ is a constant depending on $\tilde{\Delta}_0$. Consequently, $w_{\tilde{m}}=0$, $w_m >0$ for $m \in \Delta_0 \setminus \{m_o,\tilde{m}\}$. In another word, we have a toric coordinate $z$, whose components are made up with $\{\tau^{-\delta}s_m/s_{\tilde{m}}\}_{m\in S \setminus \{\tilde{m}\}}$, such that

\[
\partial \bar{\partial}\log \left(\sum_{m\in \tilde{\Delta}_0} \tau^{2\tilde{w}_m}|s_m|^2\right) = \partial \bar{\partial}\log \left(1 + \tau^{2\delta}|z|^2 + \sum_{m\in \tilde{\Delta}_0 \setminus S} \tau^{2\tilde{w}_m}|s_m/s_{\tilde{m}}|^2\right)
\]

and

\[
\tilde{s}_t = s_{m_o}+ t\left(s_{\tilde{m}} + \sum_{m\in \Delta_0 \setminus \{m_o,\tilde{m}\}} a_m s_m \right),\ \ {\rm where}\ |a_m| = \tau^{w_m}.
\]

For $X_t$, we will use the restriction of the toric \k metric on $P_\Delta$ (used in \cite{lag3}) with the \k potential 

\begin{equation}
\label{fa}
\rho_t = \tau^{-2\delta}\log \left(\sum_{m\in \tilde{\Delta}_0} \tau^{2\tilde{w}_m}|s_m|^2\right).
\end{equation}

One can observe that locally the situation can be reduced to the basic setting of section 2, with 

\[
\omega_{t,s} = \omega_t|_{X_{t,s}},\ \ \omega_t = \partial \bar{\partial} \rho_t = \sum_{k=0}^n dz_kd\bar{z}_k + O(\tau^+).
\]
\[
\tilde{p}_t(z) = \tilde{s}_t/s_{\tilde{m}} = z^m + tp(z),\ \ p(z) = 1 + \sum_{m\in \Delta_0 \setminus \{m_o,\tilde{m}\}} a_m s_m/s_{\tilde{m}}.
\]
\[
\check{p}(z) = \log p(z) = O(\tau^{\alpha'}),\ \ \alpha' = \min_{m\in \Delta_0 \setminus \{m_o,\tilde{m}\}} w_m.
\]
\[
\tilde{p}_{t,s}(z) = z^m + t e^{s\check{p}(z)}.
\]

\begin{theorem}
\label{fb}
When assumption* is satisfied and $t$ is small enough, there exists an open neighborhood $U^{(n-1)}_t$ of $\pi^{-1}(\Delta_0^{(n-1)}) \subset \partial \Delta$ such that for any generalized special Lagrangian fibre $L_{t,0}$ in $X_{t,0}$ over $U^{(n-1)}_t$, there exist a smooth family $\{L_{t,s}\}_{s\in [0,1]}$, where $L_{t,s}$ is a generalized special Lagrangian torus in $(X_{t,0}, \omega_{t,0}, \Omega_{t,s})$ that is Hamiltonian equivalent to $L_{t,0}$. When the fibre $L_{t,0}$ varies over $U^{(n-1)}_t$, $\phi_1(L_{t,1}) \subset X_{t,1} = X_t$ will form a generalized special Lagrangian fibration over $U^{(n-1)}_t$. This fibration will coincide with the fibration over $U^{\rm top}_t$ (in theorems 5.1 and 5.2 of \cite{sl1} and theorem \ref{da}) on the overlaps.
\end{theorem}
{\bf Proof:} For $\tilde{m} \in \Delta_0^{(l)}\setminus \Delta_0^{(l-1)}$ ($l\leq n-1$), our situation locally is reduced to the basic setting of section 2, with $I' = (S \setminus \{\tilde{m}\})\cap \tilde{\Delta}_0^{(l)}$ and $I'' = S \setminus (I' \cup \{\tilde{m}\})$. For a bounded constant $c_1>1$, consider the region $\tilde{U}_{\tilde{m}}$, where $|z_k| \leq c_1$ for all $k$ and $|z_k|\geq c_1^{-1}$ for $k\in I'$. Locally $X_{t,s}$ is defined by

\[
\prod_{k\in I''} z_k = -t e^{s\check{p}(z)}\prod_{k\in I'} z_k^{-m_k},
\]

where $m = m_o - \tilde{m} = (m_0,\cdots,m_n) = (m',m'')$. Therefore in $\tilde{U}_{\tilde{m}}$, we have $|z_k| \geq Ct$ for $k\in I''$. Take suitable $\alpha$ satisfying $\gamma = \alpha' + w_{m_o}\alpha >0$, we have

\[
\|\check{p}\|\xi_\alpha + (r_0^0)^{1-\alpha} \leq C\tau^\gamma,\ \ {\rm in}\ \tilde{U}_{\tilde{m}}.
\]

Let $U_{\tilde{m}}$ be the region in $\partial \Delta$ where $\tilde{U}_{\tilde{m}}\cap X_{t,0}$ is fibred over. Let $\displaystyle U^{(n-1)}_t = \bigcup_{\tilde{m} \in\Delta_0^{(n-1)}} U_{\tilde{m}}$. Then for any generalized special Lagrangian fibre $L_{t,0}$ in $X_{t,0}$ over $U^{(n-1)}_t$, when $\tau$ is small, according to the discussion in section 2, the estimates in propositions \ref{ba} and \ref{cc} together with the implicit function theorem (theorem 3.2 in \cite{sl1}) imply the existence of the family of generalized special Lagrangians $L_{t,s}$ with respect to $(X_{t,0},\omega_{t,0},\Omega_{t,s})$.\\

Since the toric \k form depends on $\tau$, the deformation (section 4) over a top dimensional face $\sigma$ of $\partial \Delta$ need some additional comment. When doing the deformation, we use the fixed \k form for $\tau$, we modify the meaning of $t$ by fixing $\tau$ and take $t\in [0,\tau^{-w_{m_o}}]$. Then we carry out the deformation as in section 4 for parameter $t$. Such toric \k metric is not of bounded geometry. It can be turned into a metric with bounded geometry in the area of discussion via a toric transformation of order $O(\tau^{|\tilde{w}|})$. It is easy to see the conclusion of theorem \ref{da} is still true for such \k metric with perhaps a slightly smaller $\beta>0$. ($\beta$ in theorem \ref{da} may need to be adjusted by $O(\frac{|\tilde{w}|}{|w_{m_o}|})$, which is small because $|w_{m_o}| \gg |\tilde{w}|$.)\\

Let $\sigma$ be a top dimensional face of $\Delta$ that is adjacent to $\tilde{m}$. In the local toric coordinate on $\tilde{U}_{\tilde{m}}$, assume that $\sigma$ corresponds to the toric devisor $\{z_0=0\}$. It is easy to see that there exists torus $L_{t,0}$ in $\tilde{U}_{\tilde{m}}$ satisfying $\log r_k$ being bounded for $1\leq k \leq n$. Such $L_{t,0}$ is of bounded geometry. Based on this, same argument as in the proof of theorem \ref{eb} will imply that $U_{\tilde{m}}$ touch $U^{\rm top}_t \cap \sigma$ in theorem \ref{da} and the two fibrations coincide on the overlap.\\
\begin{flushright} \rule{2.1mm}{2.1mm} \end{flushright}
{\bf Remark:} The result in this section can be extended to Calabi-Yau hypersurfaces in general toric varieties. The assumption* in theorem \ref{fb} can be removed with a suitable construction of the toric metric $\omega_t$ similar to the construction in \cite{toric}. Such construction is straightforward using Legendre transformation but is somewhat tedious. We will omit such construction here and defer the general case to a sequel of this paper, where using a different method, we will give a stronger construction of thin torus near $\pi^{-1}(\Delta_0^{(n-1)}\setminus \Delta_0^{(0)})$ under any fixed toric metric that does not need to depend on $t$.\\

Combining theorems \ref{da} and \ref{fb}, we have
\begin{co}
For the Calabi-Yau family $\{X_t\}$ in $P_\Delta$ satisfying assumption* with the restriction of the family of toric \k metrics $\omega_t$ from $P_\Delta$, when $t$ is small enough, there exists a smooth monodromy representing generalized special Lagrangian torus (partial) fibration for $X_t$.
\end{co}
\begin{flushright} \rule{2.1mm}{2.1mm} \end{flushright}
\begin{co}
\label{fc}
For the quintic Calabi-Yau family $\{X_t\}$ in $\mathbb{CP}^4$ that determine a primitive or a standard simplicial decomposition for $\Delta_0$, with the restriction of the family of toric \k metrics $\omega_t$ from $\mathbb{CP}^4$, when $t$ is small enough, there exists a smooth monodromy representing generalized special Lagrangian torus (partial) fibration for $X_t$.
\end{co}
\begin{flushright} \rule{2.1mm}{2.1mm} \end{flushright}

\ifx\undefined\bysame
\newcommand{\bysame}{\leavevmode\hbox to3em{\hrulefill}\,}
\fi

\end{document}